\documentclass[a4paper,12pt,reqno]{amsart}
\usepackage{amssymb}
\usepackage{amsmath}

\def\i{\,\lrcorner\,}
\def\a{\alpha}

\def\la{\langle}
\def\ra{\rangle}
\def\.{\cdot}

\def\n{\nabla}

\def\beq{\begin{equation}}
\def\eeq{\end{equation}}
\def\bea{\begin{eqnarray*}}
\def\eea{\end{eqnarray*}}

\def\r{\end{proof}}

\def\qk{quaternion--K{\"a}hler }
\def\pr{\mathrm{pr}\,}

\def \RM{\mathbb{R}}

\def \ZM{\mathbb{Z}}
\def \CM{\mathbb{C}}

\def\d{{\delta}}

\def\es{\,\lrcorner\,}
\def\we{\wedge}

\def\Ric{\mathrm{Ric}}

\def\id{\mathrm{id}}
\def\be{\begin{equation}}

\def\ee{\end{equation}}

\def\pr{{\rm pr}}
\def\tr{\mathrm{tr}}

\def\Hol{\mathrm{Hol}}
\def\hol{\mathfrak{hol}}
\def\so{\mathfrak{so}}
\def\sp{\mathfrak{sp}}
\def\R{\mathbb{R}}
\def\S{\mathrm{Sym}}

\def\End{\mathrm{End}}
\def\Cas{\mathrm{Cas}}
\def\Sp{\mathrm{Sp}}
\def\SO{\mathrm{SO}}

\def\End{\mathrm{End}}
\def\ad{\mathrm{ad}}
\def\Re{\mathrm{Re}}
\def\Aut{\mathrm{Aut}}

\def\g{{\mathfrak g}}
\def\gl{{\mathfrak gl}}

\def\we{\,{\wedge}\,}
\def\bea{\begin{eqnarray*}}
\def\eea{\end{eqnarray*}}
\def\<#1,#2>{\langle#1,#2\rangle}

\newtheorem{ede}{Definition}[section]

\newtheorem{epr}[ede]{Proposition}

\newtheorem{ath}[ede]{Theorem}

\newtheorem{elem}[ede]{Lemma}

\newtheorem{ere}[ede]{Remark}

\title{Killing forms on Quaternion--K{\"a}hler manifolds}
\author{Andrei Moroianu and Uwe Semmelmann}
\thanks{The second--named author was supported by the 
{\sl European Differential Geometry 
Endeavour} (EDGE), Research Training Network HPRN--CT--2000--00101, 
and the Centre de
Math{\'e}matiques de l'Ecole Polytechnique during the
preparation of this work.}
\address{Andrei Moroianu \\ CMAT\\ {\'E}cole Polytechnique \\ UMR 7640 du CNRS
\\ 91128 Palaiseau \\ France}
\email{am@math.polytechnique.fr}

\address{Uwe Semmelmann\\  Fachbereich Mathematik  \\ Universit{\"a}t
Hamburg \\  Bundesstr. 55 \\ 
D--20146 Hamburg, Germany}
\email{uwe.semmelmann@math.uni-hamburg.de}

\begin{document}

\begin{abstract}
We show that every Killing $p$--form on a compact \qk manifold 
has to be parallel for $p\ge2$.
\vspace{0.5cm}

\noindent
2000 {\it Mathematics Subject Classification}. Primary 53C55, 58J50
\end{abstract}

\maketitle

\section{Introduction}

In this paper we continue the study of twistor forms on compact
Riemannian manifolds with non--generic holonomy initiated in
\cite{ms} and \cite{uwe}. While the cases studied in the previous
articles concern K{\"a}hler, $G_2$-- and $Spin_7$--manifolds, we now
turn our attention to the quaternion--K{\"a}hler situation, which is the
last one in the Berger list of irreducible non--locally
symmetric Riemannian structures. 

Recall that twistor (resp. Killing) $1$--forms are duals of conformal
(resp. Killing) vector fields.
Twistor $p$--forms are natural generalizations 
of twistor $1$--forms,
defined by the property that the projection of their covariant
derivative on the Cartan product of the cotangent bundle and the
$p$--form bundle vanishes, and Killing $p$--forms have the further
property of being co--closed.

The main result of this paper is the fact that every Killing $p$--form
($p\ge 2$) on a compact \qk manifold is automatically parallel (Theorem
\ref{mainth}). The techniques used in the proof are both
representation--theoretic and analytic. We first compute some 
Casimir operators for the group $\Sp(m)\cdot \Sp(1)$ which give
explicit formulas for natural algebraic operators defined on the
exterior bundle of \qk manifolds. We then introduce natural differential 
operators (similar to $d^c$ and $\delta^c$ in K{\"a}hler geometry) 
on every \qk manifold,
compute commutator relations between them, and apply
Weitzenb{\"o}ck--type formulas in order to show that every Killing form
has to be closed, and hence parallel. 

The structure of the paper is the following. In Section 2 we recall
general facts about Killing and twistor forms, in Section 3 we
describe the decomposition of the exterior bundle of a \qk manifold
(analog to the LePage decomposition on K{\"a}hler manifolds), and in
Section 4 we introduce natural algebraic and differential 
operators on the exterior
bundle of \qk manifolds and study their behaviour with respect to this
decomposition. The next section deals with some 
representation theory, and in Section 6 we prove the main
result. Some basics on Casimir operators are explained in the Appendix.

{\it Acknowledgments.} We are grateful to Gregor Weingart for
explaining us his results on the LePage decomposition on
quaternion K{\"a}hler manifolds and his helpful comments on the 
Casimir normalization. Special thanks are due to the referee for a
thorough reading of the paper and several suggestions and
improvements.

\section{Twistor and Killing Forms on Riemannian Manifolds}

Let $( V,\,\la\cdot,\cdot\ra)$ be an $n$--dimensional Euclidean vector space. 
The tensor product  $V^*\otimes\Lambda^pV^*$ has the following
$O(n)$--invariant decomposition:
\bea\label{deco1}
V^*\otimes\Lambda^pV^*
&\cong&
\Lambda^{p-1}V^* 
\oplus \Lambda^{p+1}V^* 
\oplus 
{\mathcal T}^{p,1}V^* 
\eea
where ${\mathcal T}^{p,1}V^*$ is the intersection of the kernels of wedge
and inner product maps, which can be identified with the Cartan
product of $V^*$ and $\Lambda^pV^*$. This decomposition  immediately
translates to Riemannian manifolds $ (M^n, g)$, where we have 
\begin{equation}\label{deco}
T^*M\otimes\Lambda^pT^*M
 \cong 
\Lambda^{p-1}T^*M
\oplus
\Lambda^{p+1}T^*M
\oplus
{\mathcal T}^{p,1}M
\end{equation}
with ${\mathcal T}^{p,1}M$ denoting the vector bundle corresponding to 
the vector space 
 ${\mathcal T}^{p,1}V^* $. The covariant derivative $\nabla \psi$ 
of a $p$--form $\psi$ 
is a section of $ T^*M\otimes\Lambda^pT^*M$. Its projections onto
the summands $ \Lambda^{p+1}T^*M $ and $ \Lambda^{p-1}T^*M $
are just the differential
$d\psi$ and the co--differential $\d \psi$.
Its projection onto the third summand ${\mathcal T}^{p,1}M$ defines a
natural first order differential operator $T$, called the 
{\it twistor operator}. 
The twistor operator
$ 
T:\Gamma(\Lambda^p T^*M) \rightarrow  \Gamma({\mathcal T}^{p,1}M) 
\subset 
\Gamma(T^*M\otimes\Lambda^pT^*M)
$
is given for any vector field $X$ by the following formula
$$
[ T\psi ] (X)
 := 
[\pr_{{\mathcal T}^{p,1}M}(\nabla \psi)] ( X)
 = 
\nabla_X  \psi
 - 
\tfrac{1}{  p+1}  X \es d\psi
 + 
\tfrac{1}{ n-p+1}  X \wedge \d \psi .
$$
Note that here, and in the remaining part
of this article, we identify vectors and $1$--forms using the metric.

\begin{ede}
Differential forms in the kernel of the twistor operator are called
{\em twistor forms} (or {\em conformal Killing forms} by some
authors). 
\end{ede}

A $p$--form $ u $ is a twistor $p$--form, 
if and only if it satisfies the equation
\begin{equation}\label{killing1}
\nabla_X  u  = 
\tfrac{1}{  p+1} X\es d u   - 
\tfrac{1}{  n-p+1}  X \wedge \d  u ,
\end{equation}
for all vector fields $X$. In the physics literature this equation is
also called {\it Killing--Yano equation}.
In this article we are interested in twistor forms which are in
addition co--closed.

\begin{ede}
A $p$--form $  u  $ is called a {\it Killing $p$--form} if and only if 
$  u  $ is co--closed and in the kernel of $ T$,~i.e. 
if and only if $  u  $ satisfies 
\begin{equation}\label{killing2}
\nabla_X  u  = 
\tfrac{1}{  p+1} X\es d u  
\end{equation}
for all vector fields $X$. Clearly a Killing form is parallel if and
only if it is closed.
\end{ede}

Equivalently Killing $p$--forms may be described as $p$--forms
$ u $ for which $\nabla  u  $ is a $(p+1)$--form, or by the condition
that $X \lrcorner \nabla_X  u  = 0$ for all vector fields $X$.
Equation (\ref{killing2}) is a natural generalization of the defining equation
for Killing vector fields,~i.e. Killing $1$--forms are dual to Killing
vector fields. 

It is easy to see that $T^*T$ is an elliptic operator. Hence the space of
twistor forms is finite dimensional on compact manifolds. It actually turns
out that this space 
is finite dimensional on any connected manifold. The upper bound
of the dimension is given by the dimension of the space of twistor forms on
the standard sphere (cf.~\cite{uwe}), which coincides with the
eigenspace of the 
Laplace operator on $p$--forms corresponding to the smallest eigenvalue. In
particular, 
Killing forms on the standard sphere are precisely the co--closed minimal 
eigenforms.

The only other known examples of compact manifolds admitting Killing forms in 
degree greater than one are Sasakian manifolds, nearly K{\"a}hler manifolds,
weak $G_2$--manifolds and products of these manifolds (cf.~\cite{uwe}).

On compact manifolds one can characterize Killing  vector fields as
divergence--free 
vector fields in the kernel of $ \Delta - 2 \Ric$. A similar characterization
of arbitrary Killing forms may be given (see for instance \cite{uwe}):

\begin{epr}\label{integrab}
Let $(M^n, g)$ be a compact Riemannian manifold with a co--closed $p$--form
$ u $. Then $ u $ is a Killing form if and only if
$$
\Delta    u 
  =  
\frac{p+1}{p} q(R)    u ,
$$
where $q(R)$ is defined as the curvature term appearing in the
Weitzenb{\"o}ck formula\\ 
$\Delta = \nabla^*\nabla + q(R)$, for the Laplace operator $\Delta$
acting on $p$--forms. 
\end{epr}

In the following we need further information on the curvature term
$q(R)$. First of all it is a  symmetric
endomorphism of the bundle of differential forms defined by
\begin{equation}\label{qr}
q(R) = \sum e_j \wedge e_i \es R_{e_i,e_j},
\end{equation}
where $ \{e_i\} $ is any local orthonormal frame and $R_{e_i,e_j}$
denotes the curvature of the form bundle. On forms of degree one and
two one has an explicit expression for the action of
$q(R)$,~e.g. if $\xi$ is any $1$--form, then $q(R) \xi = \Ric(\xi)$.
In fact it is possible to define $q(R)$ in a more general
context. For this we first rewrite equation~(\ref{qr}) as
$$
q(R)
 = 
\sum_{i < j} 
(e_j \wedge  e_i \es  -  e_i \wedge  e_j\es)  R_{e_i,e_j}
 = 
\sum_{i < j} 
(e_i \we e_j)\bullet R(e_i \we e_j )\bullet
$$
where the Riemannian curvature $R$ is considered as element of
$  \S^2(\Lambda^2 T_pM) $ and $ \bullet $ denotes the standard representation
of the Lie algebra $ \so(T_pM)\cong\Lambda^2 T_pM  $ on the space of
$p$--forms. Note that we can replace $ e_i \we e_j $ by any orthonormal
basis  of $\so(T_pM)$.
Let $(M, g)$ be a Riemannian manifold with holonomy group $\Hol$. Then
the curvature tensor takes values in the Lie algebra 
$ \hol   $ of the 
holonomy group,~i.e. we can write $q(R)$ as 
$$
q(R)
 = 
\sum   \omega_i \bullet R(\omega_i)\bullet \qquad \in  \S^2(\hol)
$$
where $\{\omega_i\}$ is any orthonormal basis of $  \hol $ and $ \bullet $
denotes the form representation restricted to the holonomy group.
Writing the bundle endomorphism $q(R)$ in this way has two
immediate consequences: we see that $q(R)$ preserves any
parallel sub--bundle of the form bundle and it is clear that
by the same definition $q(R)$ gives rise to a symmetric endomorphism
on any associated vector bundle defined via a representation of
the holonomy group.

As a corollary to Proposition \ref{integrab} together with the considerations
above, we immediately obtain an
interesting property  
of Killing forms on manifolds admitting a parallel form.

\begin{elem}\label{paral}
If $\Omega$ is a parallel $k$--form and $ u $ is a Killing $p$--form on a
compact manifold $M$, then 
the contraction $ \Omega \es u$ of $ u $ with $\Omega$ is
a parallel $(p-k)$--form. 
\end{elem}
\begin{proof}
First of all we note that $\Omega \es u$ is again a Killing form. Indeed we have
$$ 
X \es \nabla_X (\Omega \es  u ) =X \es \Omega \es \nabla_X  u  = (-1)^k 
\Omega \es X \es \nabla_X  u  = 0 . 
$$
From Proposition~\ref{integrab} it follows that $\Delta  u 
=\frac{p+1}{p} q(R)  u  $. Since  
the contraction with a parallel form commutes with the Laplace operator and
with $q(R)$ we obtain 
$$
\Delta (\Omega \es  u ) 
 = 
\frac{p+1}{p} q(R) (\Omega \es  u ) .
$$ 
But since $\Omega \es  u $ 
is a Killing $(p-k)$--form, Proposition~\ref{integrab} also implies that 
$$
\Delta (\Omega \es  u ) 
 = 
\frac{p-k+1}{p-k} q(R) (\Omega \es  u ) .
$$ 
Comparing these two 
equations for $\Delta (\Omega \i  u )$ yields that the Killing form
$\Omega \es  u $ is harmonic. 
Since $M$ is compact, a harmonic form is closed, so $\Omega\i u $ is
a closed Killing form and thus parallel.
\end{proof}


\section{Exterior Forms on Quaternion--K{\"a}hler Manifolds}

 Let $(M^{4m}, g)$ be a {\it quaternion--K{\"a}hler manifold} defined by the
 holonomy reduction to $\Sp(m)\cdot \Sp(1) \subset\SO(4m)$. As usual,
 we suppose that $m\ge2$ since  $\Sp(m)\cdot \Sp(1)$ is not a proper
 subgroup of $\SO(4m)$ for $m=1$, so the holonomy condition would be
 empty in dimension $4$. 
 Any representation of~$\Sp(m)\cdot \Sp(1)$ gives rise to
 a vector bundle over $M$. The tensor product of a 
 $\Sp(1)$ representation and a $\Sp(m)$ representation defines a vector 
 bundle if and only if it factors through the projection 
 $\Sp(m) \times \Sp(1) \rightarrow \Sp(m)\cdot\Sp(1)$.
 In particular the standard representations
 $H$ resp. $E$ of $\Sp(1)$ resp. $\Sp(m)$ induce only locally defined
 vector bundles, whereas $\S^2H$ or $E \otimes H \cong T^*M^\CM$
 are globally defined. In the following we will often make no difference
 between a vector bundle and its defining representation.

More explicitly the quaternion--K{\"a}hler reduction may be described
by three locally defined almost complex structures $J_\alpha$ which
satisfy the quaternion relations and span a rank three sub--bundle
of $\End(TM)$ preserved by the Levi--Civita connection. The complexification
of this sub--bundle is isomorphic to $\S^2H$ and in any point $x\in M$
the subspace in  $\End(T_xM)$ spanned by the three almost complex
structures is isomorphic to the Lie algebra $\sp(1)$. Since
$J_\alpha$ is a skew--symmetric endomorphism one may realize 
$\sp(1)$ also as a subspace of $\Lambda^2 (T^*_xM)$. Under this
identification $J_\alpha$ is mapped to the $2$--form
$\omega_\alpha= \frac12\sum e_i \wedge J_\alpha e_i$, where
$\{e_i\}$ is any orthonormal basis of $T_xM$.

We still need some information on the decomposition of the form bundle
of a quaternion K{\"a}hler manifold. Of course this decomposition corresponds
to the decomposition of $\Lambda^p(H\otimes E)$ under the action of the
group $\Sp(1)\cdot \Sp(m)$. Details for this can be found in \cite{gregor}.
First of all it is easy to see that all possible irreducible summands
are of the form $ \S^kH \otimes \Lambda^{a,b}_0E$, where 
$\Lambda^{a,b}_0E  \subset \Lambda^{a}_0E \otimes \Lambda^{b}_0E$ denotes
the Cartan product of the two irreducible $\Sp(m)$--representations
$\Lambda^{a}_0E$ and  $\Lambda^{b}_0E$. E.g. the decomposition of the
$2$-- and $3$--forms is given as
\bea
\Lambda^2(H\otimes E) &=& 
\Lambda^{1,1}_0E  \oplus 
\S^2H \otimes [\Lambda^2_0E \oplus \CM] ,\\[1.5ex]
\Lambda^3(H\otimes E) &=&
H\otimes [\Lambda^{2,1}_0E \oplus E] \oplus 
\S^3H \otimes [\Lambda^3_0E \oplus E] .
\eea

Analyzing the form representation in more detail it is actually possible
to determine for which numbers $(k,a,b)$ the summand $ \S^kH \otimes
\Lambda^{a,b}_0E $ 
appears in the decomposition of the space of $p$--forms. In this article we
will only need the following weaker information.

\begin{elem}[\cite{gregor}]\label{forms1}
Let $ \S^kH \otimes \Lambda^{a,b}_0E $ be an irreducible summand appearing in
the decomposition of $\Lambda^p(H\otimes E)$. Then the numbers
$(k,a,b)$ satisfy 
the conditions
\bea
(i) \qquad &&        0 \le b \le a \le m ,\\[1.5ex]
(ii) \qquad &&       2b \le \min\{p-k,4m-p-k\}  ,\\[1.5ex]
(iii) \qquad &&       2a \le \min\{p+k,4m-p+k\}.
\eea
Moreover, the numbers $k,\ p$ and $a+b$ have the same parity.
\end{elem}

\bigskip

Like in the case of K\"ahler manifolds, it is possible to describe the 
action of the differential $ d $ and  
the co--differential $ \d $ on the irreducible sub--bundles $ \S^kH
\otimes \Lambda^{a,b}_0E $ of the form bundle.

\begin{elem}\label{forms2}
For any $p$--form $u$ which is a section in $ \S^kH \otimes
\Lambda^{a,b}_0E $ the forms 
$du$ and $\d u$ are sections in a sum of bundles of the type 
$ \S^{k'}H \otimes \Lambda^{a', b'}_0E $ where 
the numbers $(k',a',b')$ satisfy
the conditions
\bea
(i) \qquad && |k-k'| =  1,\\[1.5ex]
(ii) \qquad && |a - a'| + |b - b'| =1.
\eea
\end{elem}
\begin{proof}
The differential $d$ resp. the co--differential $\d u$ are projections
of $\nabla u$ onto 
$\Lambda^{p+1}T^*M$ resp. $\Lambda^{p-1}T^*M$ considered as sub--bundles
of the tensor product 
$\Lambda^{p+1}T^*M\otimes T^*M$. Hence, the conditions of the lemma follow from
the decomposition of $ (\S^kH \otimes \Lambda^{a,b}_0E) \otimes
(H\otimes E)$ into irreducible 
summands. The Clebsch--Gordan formula for $ \Sp(1)$ implies that
$\S^kH \otimes H \cong \S^{k+1}H\oplus \S^{k-1}H$, which proves the
condition on $k$.  
Similarly, the decomposition of $\Lambda^{a,b}_0E\otimes E$
implies conditions (ii) and (iii)  
(cf.~\cite{gu}).
\end{proof}


\section{Natural Operators on Quaternion-K\"ahler Manifolds}

On any K{\"a}hler manifold one can define three endomorphisms of the form
bundle: the wedge product and  
the contraction with the K{\"a}hler form and the endomorphism induced
(as a derivation) by the complex structure. 
Similarly we may associate
with any almost complex structure $J_\a$, given by the \qk reduction,
three locally defined  
endomorphisms of the form bundle. In the quaternion K{\"a}hler case the
corresponding operators  
$L_\a$, $\Lambda_\a$ and $J_\a$ are defined by
$$
L_\a:=\frac12\sum_i e_i\wedge J_\a(e_i)\wedge\qquad
\Lambda_\a:=-\frac12\sum_i e_i\es J_\a(e_i)\es\qquad
J_\a:=\sum_i J_\a( e_i)\wedge e_i\es$$
where $\{e_i\}$ is a local orthonormal base of the tangent bundle.
It is straightforward to check the following relations:
\beq\label{ini}[X\wedge,\Lambda_\a]=-J_\a(X)\es\qquad [X\wedge,J_\a]=-J_\a(X)\wedge
\eeq
\beq\label{init}[X\es,L_\a]=J_\a(X)\wedge\qquad [X\es,J_\a]=-J_\a(X)\es
\eeq

\noindent
Composing the local endomorphism $L_\alpha, \Lambda_\alpha$ and
$J_\alpha$ defined above 
we obtain globally defined endomorphisms of the form bundle: 
$$
L:=\sum_{\a} L_\a \circ L_\a,\ \qquad
L^-:=\sum_{\a} L_\a\circ  J_\a, \ \qquad 
J:=\sum_{\a} J_\a\circ  J_\a,
$$

$$
\Lambda:=\sum_{\a} \Lambda_\a\circ  \Lambda_\a,\ \qquad 
\Lambda^+:=\sum_{\a} \Lambda_\a\circ  J_\a, \ \qquad 
C:=\sum_{\a} L_\a\circ  \Lambda_\a.
$$

It is easy to prove that $J$ and $C$ are self--adjoint, while 
$L$ and $L^-$ are adjoints of $\Lambda$ and $\Lambda^+$ respectively. Moreover
it is important to 
note that $J$ and $C$ are commuting endomorphisms. The commutators of
these operators 
with the inner and wedge product with vectors are given by the
following

\begin{elem}\label{alg}
The following relations hold:
$$
\begin{array}{llll}
& [X\wedge,\Lambda]
&=&
-2\sum_{\a} \Lambda_\a\circ  J_\a(X)\es
\qquad\qquad [X\es,L]=2\sum_{\a} L_\a\circ  J_\a(X)\wedge\\[1.5ex]
& [X\wedge,L^-]
&=&
-\sum_{\a} L_\a\circ  J_\a(X)\wedge\qquad\qquad [X\es,\Lambda^+]=
\sum_{\a} \Lambda_\a\circ  J_\a(X)\es\\[1.5ex]
& [X\wedge,\Lambda^+]
&=&
-3X\es-\sum_{\a} (\Lambda_\a\circ  J_\a(X)\wedge+J_\a\circ 
 J_\a(X)\es)\\[1.5ex]
& [X\es,L^-]
&=&
\quad 3X\wedge+\sum_{\a} (J_\a\circ  J_\a(X)\wedge-
L_\a\circ  J_\a(X)\es)\\[1.5ex]
& [X\wedge,J]
&=&
-3X\wedge-2\sum_{\a} J_\a\circ  J_\a(X)\wedge\\[1.5ex]
& [X\es,J]
&=&
-3X\es-2\sum_{\a} J_\a\circ  J_\a(X)\es\\[1.5ex]
& [X\wedge,C]
&=&
\quad \sum_{\a} L_\a\circ  J_\a(X)\wedge
\qquad\qquad [X\es,C]=3X\es+\sum_{\a} \Lambda_\a\circ  J_\a(X)\wedge
\end{array}
$$
\end{elem}
These relations easily follow from (\ref{ini}) and (\ref{init}). We 
leave the necessary verifications to the reader.

We now turn our attention toward differential
operators and define natural analogues of the exterior derivative 
and co--differential on \qk manifolds. We will only introduce those of 
the operators which will be useful to the study of our particular problem.
The general theory of natural differential operators on \qk manifolds
will (hopefully) be developed in a forthcoming paper. 

Let as before $J_\a$ be a local basis of almost complex structures and define
\medskip
\bea
d^+&:=&\quad \sum_{i,\a} L_\a J_\a(e_i)\wedge\n_{e_i},\ \ 
d^-:=\sum_{i,\a} \Lambda_\a J_\a(e_i)\wedge\n_{e_i}, \ \ 
d^c:=\sum_{i,\a} J_\a J_\a(e_i)\wedge\n_{e_i},\\[1.5ex]
\d^+&:=&-\sum_{i,\a} L_\a J_\a(e_i)\es\n_{e_i},\ \ 
\d^-:=-\sum_{i,\a} \Lambda_\a J_\a(e_i)\es\n_{e_i}, \ \ 
\d^c:=-\sum_{i,\a} J_\a J_\a(e_i)\es\n_{e_i}.
\eea

\medskip

\begin{elem}\label{comfor}
The following relations hold:

\begin{displaymath}
\begin{array}{lllll} 
&[  d,\Lambda]&=\quad2\d^-\qquad\qquad &[ \d,L]&=-2d^+\\[1.5ex]
&[ d,L^-]&=-d^+\qquad\qquad &[ \d,L^-]&=-\d^+-d^c-3d\\[1.5ex]
&[ d,\Lambda^+]&=-d^-+\d^c+3\d\qquad\qquad &[ \d,\Lambda^+]&=\quad\d^-\\[1.5ex]
&[ d,J]&=-2d^c-3d\qquad\qquad &[ \d,J]&=-2\d^c-3\d\\[1.5ex]
&[ d,C]&=\quad\d^+\qquad\qquad &[ \d,C]&=-d^-+3\d
\end{array}
\end{displaymath}
\end{elem}
\begin{proof}
All algebraic operators $ L,\ \Lambda,\ L^-,\ \Lambda^+,\ J,\  $ and $ C $
appearing in the lemma  
are parallel. The result thus follows directly from Lemma \ref{alg}.
\end{proof}


\section{Representation Theoretical Results}

Fixing an arbitrary point of the manifold we may consider the bundle
endomorphisms  
$J$ and $C$ as linear maps on the space of $p$--forms
$\Lambda^p(H\otimes E)$. Obviously $J$ and $C$ are invariant under the action
of the holonomy group $\Sp(1)\cdot\Sp(m)$. Hence, both maps restrict to
equivariant maps on the irreducible sub--representations $ \S^kH\otimes
\Lambda^{a,b}_0E$. 
By the Schur Lemma, these restricted maps have to be certain multiples of
the identity. In the remaining part of this section we will show how to
compute the action of $J$ and $C$ on the irreducible components of the form
representation.

\begin{elem}\label{Jlemma}
Let $ \S^kH\otimes \Lambda^{a,b}_0 E$ be an irreducible summand of the
space of complex $p$--forms $\Lambda^p(H\otimes E)$. 
Then the bundle endomorphism $J$ acts on it as
$$
J  =  - k(k+2) \id .
$$
\end{elem}
\begin{proof}
We will consider $J=\sum J_\alpha^2 $ as a linear map acting on an
irreducible subspace  
$ \S^kH\otimes \Lambda^{a,b}_0 E $ of the space of $p$--forms
$ \Lambda^p(H\otimes E)$. 
Let $\omega_\alpha$ be the $2$--form corresponding to $J_\alpha$,~i.e. 
$\omega_\alpha= \frac12\sum e_i \wedge J_\alpha e_i$. It immediately follows
that $ J_\alpha = \omega_\alpha \bullet $ as endomorphism on the
space of $p$--forms. Here $\bullet$ denotes the
the standard representation of the Lie algebra $ \so(T_pM)\cong\Lambda^2 T_pM  $ 
on the space of $p$--forms, which is  defined as 
$$
  (X \wedge Y)  \bullet = Y \wedge X \es -X \wedge Y \es  .
$$

Let $ \g $ be a semi--simple Lie algebra equipped with an invariant 
scalar product $g$. Then the Casimir operator $ \Cas^g_\pi$, acting
on a representation $ \pi $ of $ \g$, is defined as
$  
\Cas^g_\pi = \sum \pi(X_i) \pi(X_i),   
$
where $\{X_i\}$ is an $g$--orthonormal basis of $\g$. More information
about Casimir operators can be found in the appendix.

Realizing $\sp(1)$ as a subspace of the space of $2$--forms we obtain a
scalar product on $\sp(1)$ by restricting the standard scalar product  
on $2$--forms. The corresponding Casimir operator of $\sp(1)$ is denoted by
$\Cas^{\Lambda^2}_\pi$. With respect to this standard scalar product 
the $2$--forms $\omega_\alpha$ are orthogonal and 
of length $2m$. Hence, we have for the representation
$ \pi = \S^kH\otimes \Lambda^{a,b}_0 E$ the expression
$$
J 
  =  
2m \sum_\alpha \frac{\omega_\alpha}{\sqrt{2m}}\bullet \frac{\omega_\alpha}{\sqrt{2m}}\bullet
  =  
2m    \Cas^{\Lambda^2}_\pi  .
$$
Note that $\Sp(1) \subset \Sp(1)\cdot\Sp(m)$ acts only on the
$\S^kH$--factor.

In the appendix we show the relation
$  \Cas^{\Lambda^2}_\pi = \frac{8}{2m} \Cas^{g_B}_\pi$,
where $ \Cas^{g_B}_\pi $ denotes the $\Sp(1)$--Casimir operator
defined with respect to the scalar product induced by the Killing
form. Moreover, we show that $ \Cas^{g_B}_\pi $ acts on $\S^kH$  
as $ - \frac18  k(k+2) \id$. Thus it follows that
$$
J 
  =  
2m    \Cas^{\Lambda^2}_\pi
  =  
8 \Cas^{g_B}_\pi 
  =   
- k (k+2) \id
 .
$$
We may check this formula in the case $k=1$. Here $J$ acts as
a sum of the squares of the three almost complex structures $J_\alpha$.
Hence, $J = - 3   \id $ on $ TM^\CM \cong H\otimes E$, which agrees
with our formula.
\end{proof}

\bigskip

\begin{elem}\label{Clemma}
Let $ \S^kH\otimes \Lambda^{a,b}_0 E$ be an irreducible summand of the
space of $p$--forms $\Lambda^p(T^*M^\CM)=\Lambda^p(H\otimes E)$. 
Then the bundle endomorphism $C$ acts as
$$
C =  
\frac14 
(
p(4m-p+6)  -  k(k+2)  -  4b  +  2a^2  +  2b^2  -  4(a+b)(m+1)
)
 \id  ,
$$
where $m$ is the quaternionic dimension of $M$.
\end{elem}
\begin{proof}
Starting directly from the definition of $C$ we obtain in a first step the
formula
\bea
4 C 
&=&
4 \sum   \Lambda_\alpha \circ L_\alpha
 = 
- \sum  
e_i \wedge J_\alpha (e_i)  \wedge   e_j \es J_\alpha (e_j)  \es
\\[.8ex]
&=&
\sum (e_j \wedge e_i \es  )    (J_\alpha e_j \wedge J_\alpha e_i \es  )
 + 3 \sum  e_j \wedge e_j \es
\\[.8ex]
&=&
\frac12 \sum 
(e_i \wedge e_j) \bullet (J_\alpha e_i \wedge J_\alpha e_j)\bullet 
 +  
\sum   
(e_j \wedge e_i \es  )  (J_\alpha e_i \wedge J_\alpha e_j \es  )
 + 
3p \id
\\[.8ex]
&=&
\sum_{i<j} 
(e_i \wedge e_j) \bullet (J_\alpha e_i \wedge J_\alpha e_j)\bullet 
 + J + 6p \id
 .
\eea

We now want to express the first summand in terms of the $\Sp(m)$--Casimir
operator $\Cas^{\Lambda^2}_\pi$. Hence we have to rewrite the first
summand using a basis of $\sp(m) \subset \so(4m)$. A projection map
$$
\pr: \Lambda^2(H\otimes E)  \longrightarrow  \S^2E \subset
\Lambda^2(H\otimes E) 
$$
can be defined by
$$
\pr(X \wedge Y)  =  \frac14 (X \wedge Y  +  \sum J_\alpha X \wedge
J_\alpha Y )  . 
$$
Note that  $\pr$ indeed satisfies the condition $\pr^2=\pr$. Substituting
$\pr$ into the 
formula for $4C$ we obtain
\bea
4 C 
&=&
4 \sum_{i<j}  
(e_i \wedge e_j) \bullet \pr (e_i \wedge e_j)\bullet 
 - 
\sum_{i<j}  
(e_i \wedge e_j) \bullet (e_i \wedge e_j)\bullet
 + J + 6p \id
\eea

The second summand is just the $\SO(4m)$--Casimir operator acting
on $p$--forms and a short calculation shows that it is equal to
$-p(4m-p) \id$. Moreover, we can replace the orthonormal basis 
$e_i \wedge e_j$ with any basis of  $\Lambda^2 T^*M$, which is
adapted to the decomposition of $\Lambda^2(H\otimes E)$ as
$\Sp(1)\cdot \Sp(m)$--representation and which is orthonormal
with respect to the standard scalar product of $\Lambda^2$.  
Hence, it remains only the sum over an orthonormal basis
$\{\omega_i\}$ of the summand corresponding to $\S^2E$ and 
\bea
C 
&=&
\Cas^{\Lambda^2}_{\pi}
 + 
\frac14 p(4m-p)  \id  +  \frac14 J + \frac{6}{4}p \id   ,
\eea

\noindent
where $\Cas^{\Lambda^2}_{\pi}$ denotes the $\Sp(m)$--Casimir
operator  defined with respect to the standard scalar product of
$\Lambda^2$ and acting on the representation 
$ \pi = \S^kH\otimes \Lambda^{a,b}_0 E$. Note that the action
on the $ \S^kH$--factor is trivial.


In the appendix we
show that $\Cas^{\Lambda^2}_{\pi}= 2(m+1) \Cas^{g_B}_{\pi}$, where
$ \Cas^{g_B}_{\pi} $ is the Casimir operator defined with respect
to the scalar product induced by the Killing form. Moreover we
calculate the action of $ \Cas^{g_B}_{\pi} $ on the 
$\Sp(m)$--representation $ \pi=\Lambda^{a,b}_0E $, which implies
for the Casimir operator in the $\Lambda^2$--normalization the
formula
$$
\Cas^{\Lambda^2}_{\pi}
 = 
 - \frac{1}{2} (2b-a^2-b^2+2(a+b)(m+1))  .
$$
Substituting this into our last expression for $C$ concludes the proof of
the lemma. 
\end{proof}


\section{Killing Forms on Quaternion--K\"ahler Manifolds}

The goal of this section is to prove the following 

\begin{ath}\label{mainth} Every Killing $p$--form $(p\ge2)$ on a compact \qk 
manifold $M^{4m}$ is parallel.
\end{ath}  

\begin{proof} 

Let $u$ be a Killing $p$--form on $M$. We will prove that $u$ is
closed for $p\ge 2$ and thus parallel (by (\ref{killing2})). We start
with the following

\begin{elem}\label{l62}
The exterior forms
$$d^+u,\ \d^-u,\ \d^cu,\ d^-u,\ \hbox{and}\ \d^+u+d^cu+3du$$
all vanish identically. 
\end{elem}

\begin{proof}
In relation (\ref{killing2}) we perform three operations: 1. take the
interior or wedge  product with
$J_\a(X)$ for some $\a$; 2. apply one of the operators $L_\a$, $\Lambda_\a$
or $J_\a$ to both terms; 3. sum over $\a=1,2,3$ and an orthonormal basis
$X=e_i$. This yields the following six equations: 

\begin{displaymath}
\begin{array}{lll}d^+u=\frac{1}{p+1}L^-du \qquad& d^cu=\frac{1}{p+1}Jdu\qquad
& d^-u=\frac{1}{p+1}\Lambda^+du\\[1.5ex]
\d^+u=-\frac{2}{p+1}Cdu\qquad & \d^cu=-\frac{2}{p+1}\Lambda^+du\qquad
& \d^-u=-\frac{2}{p+1}\Lambda du
\end{array}
\end{displaymath}

\bigskip

\noindent
From Lemma \ref{comfor} we then obtain

$$\begin{array}{lll}pd^+u=d(L^-u) & (p-1)d^cu=d(Ju+3u)
& pd^-u=d(\Lambda^+u)-\d^cu\\[1.5ex]
(p-1)\d^+u=-2d(Cu) & (p-1)\d^cu=-2d(\Lambda^+u)-2d^-u
& (p-3)\d^-u=-2d(\Lambda u)
\end{array}$$

As $p>1$, the third and fifth equations together show that $d^-u$ and 
$\d^cu$ are exact forms. Thus, all 6 natural first order differential 
operators $d^\pm, d^c, \d^\pm, \d^c$ map $u$ to an exact form. 
On the other hand, the right hand side equations in Lemma \ref{comfor}
show that the images of $u$ through $d^+,\d^-,\d^c,d^-$ and $\d^++d^c+3d$ 
are co--exact (in the image of $\d$). Since $M$ is compact, a form which is
simultaneously exact and co--exact must vanish.

\r

Using the second and fourth
equations above together with Lemma \ref{comfor} again, we get:
\bea
d(Ju+3u)
&=&
(p-1)d^cu  =  -\frac{p-1}{2}([d,J]u+3du),\\[1ex]
-2d(Cu) 
&=&
(p-1)\d^+u  =  (p-1)[d,C]u.
\eea

\noindent
Together with $\d^+u+d^cu+3du=0$, this gives the following system:
\beq\label{sys}
\left\{
\begin{array}{c}(p+1)d(Ju+3u) = (p-1)Jdu\\[.7ex]
(p+1)d(Cu) = (p-1)Cdu\\[.7ex]
-2Cdu+Jdu+3(p+1)du = 0
\end{array}
\right.
\eeq

The Killing form $u$ decomposes according to the decomposition of
$\Lambda^p(H\otimes E)$ 
into irreducible summands of the type $\S^kH\otimes
\Lambda^{a,b}_0E$. A priori it is not 
clear whether the irreducible components of $u$ are again Killing
forms. However, the following 
weaker statement holds:

\begin{elem}\label{cj}The forms $Ju$ and $Cu$ are again Killing forms on $M$.
\end{elem}

\begin{proof} As $d^-u=\d^cu=\d u=0$, Lemma \ref{comfor} shows that 
$Ju$ and $Cu$ are co--closed. On the other hand, $C$ and $J$ are parallel 
operators, so they commute with the curvature operator $q(R)$, 
with $\n^*\n$ and hence with the Laplace operator $\Delta$.
The characterization of Killing forms given in
Proposition~\ref{integrab} then immediately implies the result.
\end{proof}

The operators $C$ and $J$ are thus commuting 
self--adjoint linear operators acting 
on the finite dimensional space of Killing forms, so can be  
simultaneously diagonalized. We thus may assume that $u$ is an
eigenvector for both  
operators:  $Ju=ju$ and
$Cu=cu$ for some real constants $j$ and $c$.

\begin{ere}
The referee noticed, using a nice geometrical argument, that one may
actually assume that $u$ and $du$ belong to isotypical components 
$\S^kH\otimes 
\Lambda^{a,b}_0E$ and $\S^{k'}H\otimes
\Lambda^{a',b'}_0E$ of the form bundle. This information is of course
stronger than what we obtained above, but we decided, however, not to
reproduce his argument here since it is not crucial for our proof (and
Killing forms turn out to be parallel anyway!).   
\end{ere}

From now on we suppose that $du\ne 0$ and show that this leads to a
contradiction. The system (\ref{sys}) shows that $du$ is an
eigenvector of $C$ and $J$ with eigenvalues $c'$ and $j'$ satisfying
the following numerical system
\beq\label{sys1}
\left\{
\begin{array}{l}(p+1)(j+3) = (p-1)j'\\[.7ex]
(p+1)c = (p-1)c'\\[.7ex]
-2c'+j'+3(p+1) = 0
\end{array}\right.
\iff\left\{
\begin{array}{l}(p+1)(j+3) = (p-1)j'\\[.7ex]
c=\frac{j+3p}{2}\\[.7ex]
c'=\frac{j'+3(p+1)}{2}
\end{array}
\right.
\eeq

Since $u$ and $du$ do not vanish identically, there exist two 
(possibly not unique) triples of integers
$(a,b,k)$ and $(a',b',k')$
satisfying the conditions of Lemmas
\ref{forms1} and \ref{forms2} such that the $\S^kH\otimes 
\Lambda^{a,b}_0E$--component of $u$ and the $\S^{k'}H\otimes
\Lambda^{a',b'}_0E$--component of $du$ are non--zero. Using Lemmas
\ref{Clemma} and \ref{Jlemma} we obtain 
\beq \label{eq1} j=-k(k+2),\qquad j'=-k'(k'+2),\qquad c=P(k,a,b,p),\qquad 
c'=P(k',a',b',p+1),\eeq
where 
$$P(k,a,b,p):=\frac14 
(
p(4m-p+6)  -  k(k+2)  -  4b  +  2a^2  +  2b^2  -  4(a+b)(m+1)
)$$
denotes the eigenvalue of $C$ on the subspace $\S^kH\otimes 
\Lambda^{a,b}_0E$ of $\Lambda^pM$. In the remaining part of the proof
we will simply check the elementary fact that there exists no solution 
$(a,b,k,a',b',k')$ of (\ref{sys1})--(\ref{eq1}) satisfying the
compatibility conditions in Lemmas
\ref{forms1} and \ref{forms2}. 

By Lemma~\ref{forms2} (i), $k'=k\pm 1$, so $j'$ equals either
$-(k+1)(k+3)$ or $-(k-1)(k+1)$. 
From the first
equation of (\ref{sys1}) we get, in the case $k'=k+1$, 
$\frac{p+1}{p-1}(k+3)(k-1)=(k+1)(k+3)$, whose unique solution is 
$k=p$, $k'=p+1$, while
in the second case (where $k'=k-1$), the only solution is $k=1$, $k'=0$.

{\bf Case 1: $k=p$, $k'=p+1$.} From Lemma~\ref{forms1} (ii) we get
$b=b'=0$. Using the simplification $P(p,a,0,p)=\frac{1}{2}(p-a)(2m+2-p-a)$,
the last two equations of the system 
(\ref{sys1}) become:
\beq\label{syst1}
\left\{\begin{array}{l}(p-a)(2m+2-p-a)=-p(p-1)\\
(p+1-a')(2m+1-p-a')=-p(p+1)
\end{array}\right.
\eeq
From Lemma~\ref{forms2} (ii) we have $a'=a\pm1$. 
If $a'=a+1$, subtracting the 
two equations above yields $2(p-a)=2p$ hence $a=0$, so the first
equation becomes $2m+2=1$ which is impossible. 
If $a'=a-1$, subtracting again the two equations in (\ref{syst1}) yields
$$-2(2m+2-p-a)=2p,$$
so $2m+2=a$, thus contradicting Lemma~\ref{forms1} (i).

{\bf Case 2: $k=1$, $k'=0$.} In this case $j=-3$ and $j'=0$ so
(\ref{sys1}) becomes 
 \beq\label{syste}
\left\{\begin{array}{l}
P(1,a,b,p)=\frac{3(p-1)}{2}\\
P(0,a',b',p+1)=\frac{3(p+1)}{2}
\end{array}\right.
\eeq

\noindent
Subtracting these two equations yields
\beq\label{fin}
2m-p-2-2b'+2b+a'^2-a^2+b'^2-b^2+2(m+1)(a+b-a'-b')=0.
\eeq
Now, since by Lemma \ref{forms2} (ii) $(a',b')$ is one of the four
neighbors of $(a,b)$ in $\ZM^2$, we have four sub--cases:

{\bf a)} $(a',b')=(a+1,b)$. Then (\ref{fin}) gives $p=2a-3$, 
which contradicts the inequality $a\le\frac{p+k}{2}$ 
(Lemma \ref{forms1} (iii)).

{\bf b)} $(a',b')=(a-1,b)$. Then (\ref{fin}) reads $4m-p=2a-1$. 
By Lemma \ref{forms1} (i)
$a\le m$, hence $p\ge 2m+1$. From Lemma \ref{paral} we see that $\Lambda u$
(which is the contraction of $u$ with the Kraines form) has to be
parallel, so using Lemmas \ref{comfor} and \ref{l62} we get
$\Lambda(du)=d(\Lambda u)-2\d^-u=0$.
Now, it is well--known (see 
\cite{bo}) that $\Lambda$ is injective on $q$--forms for every $q\ge
2m+2$, so $du=0$ in this case.

{\bf c)} $(a',b')=(a,b+1)$. From (\ref{fin}) we get $p=2b-5$, so from
Lemma \ref{forms1} we can write 
$b\le a\le\frac{p+k}{2}=b-2$, which is impossible.

{\bf d)} $(a',b')=(a,b-1)$. Using (\ref{fin}) again we obtain
$4m-p=2b-3$ and Lemma \ref{forms1} (i) and (ii) yields $b\le
a\le\frac{4m-p+k}{2}=b-1$, a contradiction.

\r


\begin{appendix}\label{casimirapp}
\section{The computation of the Casimir eigenvalues}

Let $G$ be a compact semi--simple Lie group, with Lie algebra 
$\mathfrak g$ and 
Cartan sub--algebra ${\mathfrak t}\subset{\mathfrak g}$. Furthermore let $g$ be
any invariant scalar product on $ \g$,~e.g. $ g=g_B:= - B_\g $
where $ B_\g $ is the {\it Killing form}
defined as $B_\g(X,Y)=\tr(\ad_X\circ\ad_Y)$. For simple Lie
groups $G$ the Killing form is some
multiple of the {\it trace form} $ B_0$, which is for sub--algebras $\g \subset
\gl(n, \CM)$ defined 
as $B_0(X, Y) := \Re \tr(X\circ Y)$. For the group $G=\Sp(m)$ we have
$B_{\sp(m)} = (2m+2) B_0$.

Let $\pi:G \rightarrow \Aut(V)$ be a representation of $G$ on the
complex vector space $V$. 
If $\{X_i\}$ is a  basis of $\g$, orthonormal with respect to the
invariant scalar product $g$, the  
{\it Casimir operator} $\Cas^g_\pi\in \End(V)$ is defined as
$$ 
\Cas^g_\pi  :=  \sum_i \pi_*(X_i)\circ \pi_*(X_i),
$$
where $\pi_*:\g \rightarrow \End(V)$ denotes the differential of the
representation $\pi$. 
It is easy to see that the definition of $\Cas^g_\pi$ does not depend
from the chosen 
$g$--orthonormal basis $\{X_i\}$. Moreover, $\Cas^g_\pi$ is an endomorphism
of $V$ commuting 
with all endomorphisms of the form $\pi_*(X)$, where $X$ is any vector
in $\g$. More precisely, one defines the {\it Casimir element} $C:=\sum X_i^2$
as a vector in the universal enveloping algebra ${\mathcal U}\g$. It
then turns out that $C$ is in the center of ${\mathcal U}\g$.

If the representation $\pi$ is irreducible, the Schur Lemma implies
that $\Cas^g_\pi$ is some multiple of the 
identity. In fact it is possible to express this multiple in terms of
the highest weight of $\pi$. 

\begin{elem}\label{casimir}
Let $\lambda \in {\mathfrak t}^*$ be the highest weight of the
irreducible representation  
$\pi:G \rightarrow \Aut(V)$ and let $\rho$ be the half sum of the
positive roots of $\g$  
relative to a fixed Weyl chamber of $\mathfrak t$. The Casimir
operator is given as 
$  \Cas^{g_B}_\pi= - c_\pi \id_V $ with
$$
c_\pi  =  \|\lambda + \rho\|^2  -  \|\rho\|^2  = 
(\lambda, \lambda)  +  (\lambda, 2\rho)  , 
$$
where $(\cdot, \cdot)$ denotes the scalar product on ${\mathfrak t}^*$
induced by the Killing form $B$. 
\end{elem}

As an application we want to compute the Casimir
eigenvalues,~i.e. the scalars $c_\pi$ for several irreducible 
$\Sp(m)$--representations. With respect to the standard realization of
the Cartan algebra ${\mathfrak t} \cong \RM^m$ of $\sp(m)$, the weights 
$\lambda$ can be written as vectors
$\lambda =
(\lambda_1,\ldots,\lambda_m)=\sum \lambda_i \varepsilon_i$, 
where $\{\varepsilon_i\}$ is dual
to the standard basis in $\RM^m$. We are interested in the 
following representations of $\Sp(m)$:

$$
\begin{array}{lll}
& \pi =  \S^kE
\qquad &
\mbox{with highest weight}
\qquad
\lambda  =  k \varepsilon_1   =  (k,0,\ldots,0)\\[2ex]
& \pi  =  \Lambda^a_0E
\qquad &
\mbox{with highest weight}
\qquad
\lambda  =  \sum^a_{i=1} \varepsilon_i   =  (
\underbrace{1,\ldots,1}_{a},0,\ldots,0)\\[1.5ex]
& \pi  = \Lambda^{a,b}_0E
\qquad &
\mbox{with highest weight}
\qquad
\lambda  =  \sum^a_{i=1} \varepsilon_i  + \sum^b_{i=1} 
\varepsilon_i \\[1.5ex]
&&
\phantom{XXXXXXXXXXXX}   =  (
\underbrace{2,\ldots,2}_{b},\underbrace{1,\ldots,1}_{a-b},0,\ldots,0)
\end{array}
$$
Under the identification ${\mathfrak t} \cong \RM^m$, the trace 
form $B_0$ corresponds
to twice the standard scalar product on $\RM^m$. Hence, in the 
formula of Lemma~\ref{casimir}
we can replace $(\cdot, \cdot)$ by $\frac{1}{4(m+1)}$ times the 
standard scalar product on $\RM^m$. Moreover the half--sum of
positive roots is given as the vector
$ \rho = (m, m-1,\ldots,1)$.
Using these remarks we obtain the following Casimir eigenvalues:

$$
\begin{array}{lll}
& \pi =  \S^kE
\qquad & 
c_\pi  = \frac{1}{4(m+1)}  k(k+2m)\\[1.5ex]
& \pi  =  \Lambda^a_0E
\qquad &
c_\pi  =  \frac{1}{4(m+1)} a(2-a+2m)\\[1.5ex]
& \pi  = \Lambda^{a,b}_0E
\qquad &
c_\pi  =  \frac{1}{4(m+1)} (2b-a^2-b^2+2(a+b)(m+1))
\end{array}
$$

In particular we obtain for $m=1$ the Casimir eigenvalues for 
$\Sp(1)$,~e.g. the
Casimir operator on $\S^kH$ is given as $-\frac18k(k+2)\id$. 
Moreover, the Casimir eigenvalue
for $\pi = \Lambda^a_0E $ is of course the special case 
$b=0$ for $\pi =\Lambda^{a,b}_0E$.
Note that the Casimir eigenvalue (with respect to the Killing 
form) of the adjoint representation
is always one. Using our formula we can check this for 
$\sp(m)\cong\S^2E$.

In the remaining part of this section we want to make some 
remarks concerning the 
Casimir normalization,~i.e. we will give a formula comparing 
the Casimir operators
corresponding to different scalar products.

Let $\g$ be the Lie algebra of a compact simple Lie group and let $V$
be a real vector space with a isotypical representation of $\g$,
i.e. $V$ is a sum of isomorphic irreducible $\g$--representations.
Thus, the Casimir operator $\Cas^g_V$ acts on $V$ as 
$\Cas^g_V= - c^g_V \id$ for some number $c^g_V$. 
Moreover we assume
that $V$ is equipped with a $\g$--invariant scalar product 
$\la \cdot,\cdot\ra$,~i.e.
$\g \subset \so(V) \cong \Lambda^2V$. Restricting the 
induced scalar product
onto  $\Lambda^2V$ we obtain a natural scalar product 
$\la \cdot,\cdot\ra_{\Lambda^2}$ on $\g$.
Note that,
$$
\la \alpha,  \beta \ra_{\Lambda^2}
 = 
-\frac12 \tr_V(\alpha \circ \beta)
 = 
\frac12 \la\alpha, \beta\ra_{\End V}  .
$$

\begin{elem}
Let $\Cas^{\Lambda^2}_\pi$ be the $\g$--Casimir operator defined 
with respect to the scalar product  $\la \cdot,\cdot\ra_{\Lambda^2}$ 
restricted to $\g$ and let $ \Cas^g_\pi$ be the $\g$--Casimir operator
corresponding to any other invariant scalar product $g$. Then for
any irreducible $\g$--representation $\pi$ it follows that
$$
\Cas^{\Lambda^2}_\pi
 = 
2 \frac{\dim \g}{\dim V}   \frac{1}{c^g_V}     \Cas^g_\pi  .
$$
\end{elem}
\begin{proof}
Let $\{X_i\}$ be an orthonormal basis of $\g$ with respect to
$\la \cdot,\cdot\ra_{\Lambda^2}$ and let $\{e_i\}$ be an
orthonormal basis of $V$. Then
$  
\Cas^{\Lambda^2}_V (v) = - c^{\Lambda^2}_V v = \sum_i X^2_i(v)  
$
and we get

$$
-  \dim V  c^{\Lambda^2}_V 
 = 
\sum_{i,j}\la X_i^2(e_j), e_j  \ra
=
-\sum_{i,j}\la X_i(e_j), X_i(e_j)  \ra
=
-2\sum  |X_i|^2_{\Lambda^2}
=
-2\dim \g
$$
which proves the lemma in the case $\pi = V$. Since $\g$ is
a simple Lie algebra it follows that two Casimir operators
defined with respect to different scalar products differ only
by a factor independent from the irreducible representation $\pi$. 
Hence
$$
\frac{c^{\Lambda^2}_\pi}{c^g_\pi}
  =  
\frac{c^{\Lambda^2}_V}{c^g_V}
$$
and the statement of the lemma follows from the special case $\pi=V$.
\end{proof}

As a first application we consider the case $  \g=\sp(m)  $ with 
$  V=\RM^{4m}\cong E  $
and $  \dim \sp(m)= m(2m+1)$. Since $  V \cong E = \S^1 E  $
the formulas above imply that 
$  c^{g_B}_V=\frac{2m+1}{4(m+1)}$. 
Hence,

$$
\Cas^{\Lambda^2}_\pi
 = 
2 \frac{m(2m+1)}{4m} \frac{4(m+1)}{2m+1}     c^{g_B}_\pi 
 = 
2 (m+1)  \Cas^{g_B}_\pi   .
$$

As a second application we want to derive a similar formula for the
$\Sp(1)$--Casimir operators. Here we take $V$ to be the real subspace of
$H\otimes E$. Hence $V \cong \R^{4m}$ with the standard representation
of $\Sp(1)$, acting trivially on the $E$--factor. The general formula
then implies
$$
\Cas^{\Lambda^2}_\pi
  =  
2 \frac{3}{4m} \frac{8}{3}   c^{g_B}_\pi
  =  
\frac{8}{2m} \Cas^{g_B}_\pi  .
$$

\end{appendix}


 \labelsep .5cm

\end{document}


\begin{center}                 

{\bf\large Killing forms on Quaternion-K{\"a}hler manifolds}

\vs

{\bf-- Corrigendum --}

\vs

{Andrei Moroianu and Uwe Semmelmann}
\end{center}

\vs\vs

The aim of this note is to fill a gap in the proof of Theorem 6.1 in
\cite{ms} - stating that every Killing $p$-form $(p\ge2)$ on a compact \qk 
manifold $M^{4m}$ ($m\ge 2$) is parallel. This gap, pointed out by
Liana David, is due to two 
wrong coefficients in the formulas at the middle of page 329. The
correct equations read
$$\begin{array}{lll}pd^+u=d(L^-u), & (p-1)d^cu=d(Ju+3u),
& pd^-u=d(\Lambda^+u)-\d^cu,\\[1.5ex]
(p-1)\d^+u=-2d(Cu), & (p-1)\d^cu=-2d(\Lambda^+u)-2d^-u,
& (p-3)\d^-u=-2d(\Lambda u).
\end{array}$$
\ni The remaining part of the proof works verbatim for $p>3$, but an
extra argument is needed for $p=2$ and $p=3$.
\vs 

{\bf The case  $p=2$.} From the 6 equations above one obtains the
vanishing of $d^+u,\ \d^-u,$ and 
$\d^+u+d^cu+3du$, but no longer that of $\d^cu$ and
$d^-u$. Correspondingly, the proof of Lemma 6.3 fails. Fortunately, an
ad hoc argument shows that $du=0$ in this case. Indeed, the third
equation of the system $(7)$ shows that $du$ is an eigenform of $2C-J$
for the eigenvalue 9. On the other hand, Lemma 5.1, Lemma 5.2 and the
decomposition 
$$\Lambda^3(H\otimes E) =
H\otimes [\Lambda^{2,1}_0E \oplus E] \oplus 
\S^3H \otimes [\Lambda^3_0E \oplus E] $$
\ni show that the eigenvalues of $2C-J$ on the four summands of
$\Lambda^3M$ are respectively $3$, $4m+5$, $15$, and $4m+11$. For
$m\ge 2$, none of them equals 9.
\vs

{\bf The case  $p=3$.} The proof works well in this case, except that
one does not obtain $\d^-u=0$ in Lemma 6.2. However, this relation is
not needed 
until the point (b) at the bottom of page 331. In order to rule out
that case,  one has to
replace the argument given there with the fact that for $p=3$ and $m\ge
2$, the inequality $p\ge 2m+1$ is impossible.

{\bf Remark.} The assumption that $m\ge 2$ is essential. For $m=1$,
the quaternionic projective space $\HM {\rm P}^1\simeq S^4$ carries
non-parallel Killing 2-forms and 3-forms (cf. \cite{u}).


 \labelsep .5cm